\documentclass{article}
\usepackage[applemac]{inputenc}
\usepackage{amsmath,amssymb,fullpage}
\usepackage[T1]{fontenc}
\usepackage{amssymb}
\usepackage{amstext}
\input amssym.def
\input amssym.tex
\usepackage{epsfig}
\usepackage{amsmath}
\usepackage{amssymb}
\usepackage{amstext}
\usepackage{yfonts}
\usepackage{epsfig}
\usepackage[]{graphicx}
\usepackage{color}
\usepackage{comment}

\makeatletter
\DeclareFontEncoding{LS1}{}{}
\DeclareFontEncoding{LS2}{}{\noaccents@}
\DeclareFontSubstitution{LS1}{stix}{m}{n}
\DeclareFontSubstitution{LS2}{stix}{m}{n}
\makeatother

\DeclareMathAlphabet\mathscr{LS1}{stixscr}{m}{n}
\SetMathAlphabet\mathscr{bold}{LS1}{stixscr}{b}{n}
\DeclareMathAlphabet\mathcal{LS2}{stixcal}{m}{n}
\SetMathAlphabet\mathcal{bold}{LS2}{stixcal}{b}{n}

\newtheorem{example}{Example}[section]
\newtheorem{theorem}[example]{Theorem}

\newtheorem{definition}[example]{Definition}
\newtheorem{lemma}[example]{Lemma}

\newtheorem{remark}[example]{Remark}

\author{\bf Jasmina \DJ or\dj evi\'{c}, \thanks{  {\it Adress:} Faculty of Science and Mathematics, University of Ni\v s, Vi\v segradska 33, 18000 Ni\v s, Serbia, \newline The Faculty of Mathematics and Natural Sciences, University of Oslo, Blindern 0316 Oslo, Norway,      {\it E-mail:} jasmindj@math.uio.no, djordjevichristina@gmail.com, jasmina.djordjevic@pmf.edu.rs.}\,\, Andrey Dorogovtsev
\thanks{ {\it Adress:} Institute of Mathematics National Academy of Sciences of Ukraine,    {\it E-mail: andrey.dorogovtsev@gmail.com.}  }
}
\date{}
\title{Backward Stochastic Differential Equations with interaction}

\begin{document}
\maketitle

\begin{abstract}

\end{abstract}  In this paper  backward stochastic differential equations  with interaction (shorter BSDEs with interaction) are introduced. Far to our knowledge, this type of equation is not seen in the literature before. Existence and uniqueness result for BSDE with interaction is proved under version of Lipschitz condition with respect to Wasserstein 
distance.  Such kind of BSDE arises naturally when considering Monge-Kantorovich problem. In the proof we start from discrete measures using known result of Pardoux and Peng and approximate general measure via Wasserstein 
distance.

\bigskip



\section{Introduction}

In the article we consider backward stochastic differential equations with interaction. To explain our interest to this equation let us recall some facts about the equations with interaction  and mass transportation problem. Equations with interaction were introduced in \cite{andadd2}, \cite{andadd1},  \cite{and2},  \cite{anddod1},  \cite{andadd3}  in the following form
\begin{eqnarray}
 &&dx(u,t)=a(x(u,t), \mu_t,t)dt+ \int_{{\mathbb R}^d}b(x(u,t), \mu_t, t, p)W(dp, dt),\nonumber\\
&& x(u,0)=u, u\in{\mathbb R}^d,\nonumber\\
&& \mu_t=\mu_0\circ x(\cdot, t)^{-1}.\label{add1}
  \end{eqnarray}

 Here $W$ is a $d$-dimensional Wiener sheet on ${\mathbb R}^d\times [0; +\infty).$ $x(u, \cdot)$ is a trajectory in ${\mathbb R}^d$ of the particle which starts from the point $u\in{\mathbb R}^d$, $\mu_0$ is the initial mass distribution which is supposed to be a probability measure on the Borel $\sigma$-field in ${\mathbb R}^d.$ Correspondingly $\mu_t=\mu_0\circ x(\cdot, t)^{-1}$ is a mass distribution of all particles at time $t.$ In the Eq. (\ref{add1}) coefficients of drift and diffusion (interaction with the random media) depend on the measure $\mu_t.$ It means, that the infinitesimal  increments of trajectory of one particle depend not only on its position but also on the mass distribution of all system. Note, that such kind of equation is different from McKean-Vlasov equation because together with the mass distribution $\mu_t, t\geq0$ it allows to determ a flow $x(\cdot, t), t\geq0$  which carries this mass and study its properties. The detailed investigation of equations and related   Markov processes was done in the monograph \cite{and1}  
 and the articles \cite{moi2}, \cite{andadd4} etc. The main feature of Eq. (\ref{add1}), which is important for us in this paper, is the existence of the flow. This allows to connect Eq. (\ref{add1}) to Monge-Kantorovich optimal transport problem. Together with the initial mass distribuion $\mu_0$ consider the probability measure $\nu$ on ${\mathbb R}^d.$ If one treates $\mu_t$ from Eq. (\ref{add1}) as a mass distribution at the moment $t\geq0,$ then the natural question arise. Is it possible to have an equality 
\begin{equation}
 \mu_T=\nu \label{add2}
 \end{equation}
 for fixed moment $T?$ For sufficiently good measures $\mu_0$ and $\nu$ this question can be reformulated. It is known \cite{ced} , that if $\mu_0$ and $\nu$ has a density with respect to Lebesgue measure,  then there exists Brenier map $F: {\mathbb R}^d\to{\mathbb R}^d$ such that
 $$
 \nu=\mu_0\circ F^{-1}.
 $$
 This map exists under less restrictive conditions on $\mu_0$ and $\nu$ but even for smooth measures $F$  can be  nondifferentiable. This is one of the reasons why one can not expect  (\ref{add2}) for the smooth coefficients $a$ and $b.$ Due to described problem we consider in this article the following backward stochastic differential equation with interaction
\begin{eqnarray}
&& dY(u,t)=a(Y(u,t), \mu_t,t)dt+ b(Y(u,t), \mu_t,t))dW(t)+Z(u,t)dW(t),\nonumber\\
 &&Y(u, T)=\xi(u),  u\in{\mathbb R}^d,\nonumber\\
 &&\mu_t=\mu_0\circ Y(\cdot, t)^{-1}.\label{add3}
  \end{eqnarray}

In this equation random field $\xi(u), u\in{\mathbb R}^d$  plays the role of Brenier map and the additional process $Z$ arise in order to help the state process $Y$ to reach value $\xi$ at the time $T.$ The main result of the paper is the theorem about the uniqueness and existence of solutions to Eq. (\ref{add3}).

\medskip

Let us recall briefly some known facts about backward stochastic differential equations (shorter BSDEs)

BSDEs  occurred first time in 70s when  Bismuth (1973), Kushner (1972) and later on  Bensoussan (1982) and  Haussmann (1986) (see  \cite{ben}, \cite{bis},  \cite{hau}, \cite{kus}) observed stochastic optimal control problem. In order to optimize  cost function, they defined Hamiltonian and trough maximal conditions, linear form of BSDE raised.   Later on Pardoux and Peng \cite{pp}, as well as Mao \cite{mao}, generalized the BSDE and gave the existence and uniqueness result for the nonlinear BSDE (which has general drift function and additional function dependent of state process in addition to control process in diffusion coefficient).  Reader can find more results from the authors concerning stability of BSDEs in  \cite{moiadd2} -- \cite{moiadd1}, \cite{moi1}, \cite{moiadd3}  etc. We will need in this paper multi dimensional BSDEs which is analyzed by  Pardoux and Peng in \cite{pp}. In order to get multidimensional solution to BSDE 
\begin{equation}
y(t)=y(T)-\int_t^T\!\!f(s,y(s),z(s))ds-\int_t^T\!\big[g(s,y(s), z(s))+z(s)\big]\,dB(s),\label{1}
\ \ t\in[0,T],
\end{equation}
 with a terminal condition  $y(T)\in L^2_{{\cal
 F}_T}(\Omega;R^d)$ Pardoux and Peng used method of Picards iterations.  The sequence is defined in following way. Set $y_0(t)=z_0(t)\equiv 0, t\in[0,T]$, and for every $n\geq 1$  define
 \begin{equation}
y_{n+1}(t)=\xi-\int_t^T\!f(y_{n}(s),z_{n}(s),s)\,ds-\int_t^T\![g(y_{n}(s),s)+z_{n+1}(s)]\,dB(s), \ \ t\in [0,T].\label{pikarova}
\end{equation}
Eq. (\ref{pikarova}) has a unique solution  which is a pair of adapted  processes $(y_{n+1}(\cdot), z_{n+1}(\cdot))$. They further proved that this  sequence of processes   converges to the pair of process  $(y(\cdot), z(\cdot))$  which is a solution of (\ref{1}). 
 
 \medskip

Paper is organized as follows. In Section 2 preliminaries and the definition of the main problem are given. Section 3 contains proofs of main results given in three main steps. 

\section{Preliminaries and formulation of the problem}

Let  $(\Omega, \mathcal{F}, P)$  be probability space on which a standard $d$-dimensional Brownian motion $B=\left(B_{t}\right)_{0 \leq t \leq  T}$ with natural  augmented filtration $\left\{\mathcal{F}_{t}^{B}\right\}_{0 \leq t \leq T}$ is defined. Also let

\medskip


$ \mathcal{M}^2(\mathbb{R}^d\times [0,T]; \mathbb{R}^d) ( \mathcal{M}^2(\mathbb{R}^d\times [0,T]; \mathbb{R}^{d\times d}))$  be a set of $\mathbb{R}^d$-valued ($\mathbb{R}^{d\times d}$-valued) processes which are $
\mathcal{F}_t$-progressively measurable, and the square integrable with respect to $P\times \lambda$. Euclidean norm in $R^d, d\in \mathbb{N}$ is denoted by $|\cdot|$, while the trace norm of a
matrix $A$ is denoted by $||A||=\sqrt{\mbox{{\it trace}}[A^TA]}$,
where $A^T$ is the transpose of a matrix or vector. Let $\mathcal{B}^d$, $d\geq 1$  ($\mathcal{B}^{d\times d}$, $d\geq 1$) denotes the Borel sets on $\mathbb{R}^d$, $d\geq 1$,  ($\mathbb{R}^{d\times d}, d\geq 1$) and $\textfrak{M}$ be the set of all probability measures on Borel $\sigma$-field, $\mathcal{B}^d$. We introduce basic BSDE  with interaction  which is a special case of Eq.  (\ref{add3}) (as function $g(\cdot, \cdot, \cdot)\equiv 0$).

 \begin{align}\label{1in}
\left\{\begin{array}{lr}
 dY(u,t)=f(Y(u,t), \mu_t, t )dt+ Z(u,t) dB(t),\ \ t\in [0,T]\\
Y(u,T):=\xi(u),\  \  u\in \mathbb{R}^n, \quad \xi(u)\in L^2(\Omega),\\
\mu_t=\mu_0\circ Y^{-1}(\cdot,t),
\end{array}\right.
\end{align}
where $\mu_t=\mu_0\circ Y^{-1}(\cdot,t)$ is an image of $\mu_0$ under the mapping $Y(\cdot,t):\mathbb{R}^ d\rightarrow \mathbb{R}^ d$.


 

 \medskip

Let us suppose:

 \medskip

{\it {\bf (A1)}} $f: \mathbb{R}^ n\times \textfrak{M} \times [0,T] \rightarrow \mathbb{R}^d$ is  $ \mathcal{B}^n\times\mathcal{B}(\textfrak{M})\times  \mathcal{B}([0,T])$-measurable.


\bigskip

\begin{definition}\label{d1} Pair of random fields $(Y(u,t), Z(u,t))_{0\leq t \leq T}\in \mathcal{M}^2( \mathbb{R}^d\times [0,T]; \mathbb{R}^d)\times \mathcal{M}^2( \mathbb{R}^d\times [0,T]; \mathbb{R}^{d\times d})  $  is a solution to Eq.(\ref{1in}) if pair $(Y(u,t), Z(u,t))$ satisfies  integral Eq.(\ref{1in})  almost surely.

The solution $(Y(u,t), Z(u,t))_{0\leq t \leq T}\in \mathcal{M}^2( \mathbb{R}^d\times [0,T]; \mathbb{R}^d)\times \mathcal{M}^2(  \mathbb{R}^d\times [0,T]; \mathbb{R}^{d\times d})$ is called unique if for any other solution of Eq.(\ref{1in}) $(\bar{Y}(u,t), \bar{Z}(u,t))_{0\leq t \leq T}\in \mathcal{M}^2(  \mathbb{R}^d\times [0,T]; \mathbb{R}^d)\times \mathcal{M}^2( \mathbb{R}^d\times  [0,T]; \mathbb{R}^{d\times d}) $  

$$E\int_0^T\left(|Y(u,t)-\bar{Y}(u,t)|^2+||Z(u,t)-\bar{Z}(u,t)||^2\right)dt =0.$$
\end{definition}

\medskip

Let us restrict ourselves to the case when the driver, i.e. function $f$ of general BSDE with interaction (\ref{1in}) has a special form

\begin{equation}
f(y,\mu_t)= \int_{\mathbb{R}^d} h( y,u) \mu_t(du),\label{driverf}
\end{equation}
where $h:\mathbb{R}^d\times\mathbb{R}^d\rightarrow \mathbb{R}^d$.  In order to prove the existence and uniqueness of the solution to Eq. (\ref{1in}), we first suppose that measure $\mu_0$ is discrete 
 
$$\mu_0=\sum_{k=1}^N p_k \delta_{u_k},$$
for $p_k>0, k=\overline{1,N}, \sum_{k=1}^ N p_j= 1$,  $u_k\in \mathbb{R}^ d, k=\overline{1,N}.$ In this case measure $\mu_t$  has the form

$$\mu_t=\sum_{k=1}^N p_k \delta_{Y(u_k,t)},$$
and one can get BSDE  with interaction for heavy particles

\begin{equation}
 d Y(u_k,t)=\sum_{j=1}^N h(Y(u_k,t),Y(u_j,t)) p_j dt+ Z(u_k,t) dB(t),\ \   Y(u_k,T)=\xi(u_k), k=\overline{1,N},\label{heavy1}
 \end{equation}



Our analysis of introduced problem consist of two main steps.

\smallskip

\underline{\it {Step 1.}} {\bf Proof of the existence of the solution to Eq. (\ref{heavy1}),  BSDE with interaction with discrete measure for heavy particles.}

\medskip

Introduce following assumptions for function $h$ and terminal condition  $\xi$:\\

{\it {\bf (A2)} }  For any  $y, y_1,v,v_1 \in \mathbb{R}^d$

$$|h(y,v)-h(y_1,v_1)| \leq L_1(|y-y_1|+|v-v_1|),$$

for some positive constant $L_1$.

{\it {\bf (A3)} } For any  $u_1, u_2\in \mathbb{R}^d$

$$E|\xi(u_1)-\xi(u_2)|^2\leq L_2 |u_1-u_2|^2, $$

for some positive constant $L_2$. 

\noindent Without loss of generality, we will use $L=\max\{L_1,L_2\}$. We will omit notation of  dependence from time $t$ in functions in the sequel.

\smallskip

\begin{theorem}\label{t2.2}
Suppose that driver $f$ satisfies assumption $(A1)$, and function $h$ satisfies assumption  $(A2)$, then there exists a unique solution to Eq. (\ref{heavy1}).
\end{theorem}

\begin{dproof} The steps which we use are similar to the ones from the paper of Pardoux and Peng \cite{pp}. First we will prove {\it uniqueness} of the solution. Let us suppose that there exist two solutions $(Y(u_k,t), Z(u_k,t))_{0\leq t \leq T}$, \\
$(\overline{Y}(u_k,t),\overline{Z}(u_k,t))_{0\leq t \leq T}\in \mathcal{M}^2( \mathbb{R}^d\times [0,T]; \mathbb{R}^d)\times \mathcal{M}^2(  \mathbb{R}^d\times [0,T]; \mathbb{R}^{d\times d}) $ which satisfy equation (\ref{heavy1}). If we apply It\^o formula to $|Y(u_k,t)-\overline{Y}(u_k,t)|^2$ and take expectation  we obtain
\begin{eqnarray}
&&\hskip-0.5cm E|Y(u_k,t)-\overline{Y}(u_k,t)|^2+E\int_0^T||Z(u_k,t)-\overline{Z}(u_k,t)||^2dt\nonumber\\
&&\hskip-0.5cm  \qquad =-2E\int_0^T\left(\sum_{j=1}^N[ h(Y(u_k,t),Y(u_j,t)) -h(\overline{Y}(u_k,t),\overline{Y}(u_j,t))]p_j, Y(u_k,t)-\overline{Y}(u_k,t)]\right)dt\nonumber\\
&&\hskip-0.5cm \qquad \leq 2^{N-1}T E\int_0^T\sum_{j=1}^N\left(h(Y(u_k,t),Y(u_j,t)) -h(\overline{Y}(u_k,t),\overline{Y}(u_j,t))p_j\right)^2dt\\
&&\hskip-0.5cm \qquad \quad +E\int_0^T|Y(u_k,t)-\overline{Y}(u_k,t)|^2 dt \nonumber\\
&&\hskip-0.5cm  \qquad \leq C_1 \int_0^TE|Y(u_k,t)-\overline{Y}(u_k,t)|^2dt\nonumber\\
&&\hskip-0.5cm  \qquad \quad+C_2 \sum_{j=1}^N\left[\int_0^TE|Y(u_j,t)-\overline{Y}(u_j,t)|^2dt +  \int_0^TE|Y(u_k,t)-\overline{Y}(u_k,t)|^2dt\right],\label{un1}
\end{eqnarray}
where $C_1, C_2$ are constants ($C_2$ is dependent from Lipschitz constant $L$). From last expression, we have 

\begin{eqnarray*}
&&E|Y(u_k,t)-\overline{Y}(u_k,t)|^2\\
&&\qquad\qquad \leq C_1 \int_0^T E|Y(u_k,t)-\overline{Y}(u_k,t)|^2dt\nonumber\\
&&\qquad \qquad \quad+C_2 \sum_{j=1}^N\left[\int_0^T E|Y(u_j,t)-\overline{Y}(u_j,t)|^2dt +  \int_0^T E|Y(u_k,t)-\overline{Y}(u_k,t)|^2dt\right],\\
&&\qquad \qquad =(C_1+N) \int_0^T E|Y(u_k,t)-\overline{Y}(u_k,t)|^2+C_2 \sum_{j=1}^N \int_0^T E |Y(u_j,t)-\overline{Y}(u_j,t)|^2dt ,
\end{eqnarray*}
and  if sum up last inequality for all  $k=\overline{1,N}$, we obtain
\begin{eqnarray*}
&&\sum_{k=1}^N E|Y(u_k,t)-\overline{Y}(u_k,t)|^2 \leq( C_1+N)\sum_{k=1}^N \int_0^T E|Y(u_k,t)-\overline{Y}(u_k,t)|^2dt\nonumber\\
&&\qquad \qquad \qquad \qquad \qquad \qquad \qquad+C_2 \sum_{j=1}^N \int_0^T E |Y(u_j,t)-\overline{Y}(u_j,t)|^2dt.
\end{eqnarray*}
If we define 
$$a(t):=\sum_{k=1}^N E|Y(u_k,t)-\overline{Y}(u_k,t)|^2, \ t\in[0,T]$$ 
then
$$
 a(t) \leq c \int_0^Ta(s)ds,$$
for generic constant $c$. Applying  Gronwall lemma we conclude that $a(t)=0, t\in[0,T]$, and 
\begin{equation}Y(u_k,t)=\bar{Y}(u_k,t), \ t\in[0,T] \ \ a.s.\ .\label{uny}
 \end{equation}
If we substitute (\ref{uny}) in (\ref{un1}), then
$$E\int_0^T||Z(u_k,t)-\overline{Z}(u_k,t)||^2dt=0,\ \ a.s., $$
which completes the proof of uniqueness.

\medskip

Next we will prove {\it existence} of the solution for the BSDE (\ref{heavy1}) with interaction for heavy particles using Picard iterations. For every $k=\overline{1,N}$ let $Y^0(u_k,t)\equiv 0$, and $\{(Y^i(u_k,t), Z^i(u_k,t)), t\in[0,T]\}_{i\geq1}$ is a sequence in $\mathcal{M}^2(  [0,T]; \mathbb{R}^d)\times \mathcal{M}^2(   [0,T]; \mathbb{R}^{d\times d}) $ defined recursively in the following way

\begin{eqnarray*}
&&Y^{i+1}(u_k,t)= u_k+\int_t^T\sum_{j=1}^N h(Y^{i}(u_k,s),Y^{i}(u_j,s)) p_j ds+ \int_0^TZ^{i+1}(u_k,s) dB(s),\\
&&  Y^{i+1}(u_k,t)=\xi(u_k).
\end{eqnarray*}
This equation has drift in step $i+1$ which is only dependent from $t$, i.e. it is of a form

$$y(t)= u_k+\int_t^T l(s)ds+ \int_t^Tz(s) dB(s),\ \  y(T)=\xi(u_k),$$
for any $k$. This equation has unique solution by martingale representation theorem (for more details see \cite{pp} by Pardoux and Peng). Similarly as in the case of the uniqueness, we obtain

\begin{eqnarray*}
&&E|Y^{i+1}(u_k,t)-Y^{i}(u_k,t)|^2+\int_t^T||Z^{i+1}(u_k,s)-Z^{i}(u_k,s)||^2ds\\
&&\qquad \leq C_3\sum_{j=1}^N\left(\int_t^T E |Y^{i}(u_j,s)-Y^{i-1}(u_j,s)|^2ds+\int_t^TE |Y^{i}(u_k,s)-Y^{i-1}(u_k,s)|^2ds\right),
\end{eqnarray*}
for a generic constant $C_3$. If we sum up last expression for every $k$, where $k=\overline{1,N}$, we obtain

\begin{eqnarray*}
&&\sum_{k=1}^NE|Y^{i+1}(u_k,t)-Y^{i}(u_k,t)|^2 +\sum_{k=1}^N\int_t^T||Z^{i+1}(u_k,s)-Z^{i}(u_k,s)||^2ds\\
&&\hskip 5cm\leq C_4 \sum_{j=1}^N \int_t^T E|Y^{i+1}(u_j,s)-Y^{i}(u_j,s)|^2ds,
\end{eqnarray*}
for a generic constant $C_{4}$. For 
\begin{eqnarray*}
&& a_i(t):=\sum_{k=1}^N \int_t^TE|Y^{i}(u_k,t)-Y^{i-1}(u_k,t)|^2dt,\\
&&b(t):=\sum_{k=1}^N\int_t^T||Z^{i+1}(u_k,s)-Z^{i}(u_k,s)||^2ds,
\end{eqnarray*}
we have 

\begin{equation}
-\frac{e^{Kt}a_{i+1}(t)}{dt}+e^{Kt} b_{i+1}(t)\leq 0,\   \   Y_{0}(t)=0. \label{ex1}
\end{equation}
Integrating last inequality
$$a_{i+1}(t)    +   \int_t^Te^{K(s-t)}b_{i+1}(t)dt\leq 0.$$
Hence $\{(Y^i(u_k,t), Z^i(u_k,t)), t\in[0,T]\}_{i\geq1}$ is a Cauchy sequence in $\mathcal{M}^2(  [0,T]; \mathbb{R}^d)\times \mathcal{M}^2(   [0,T]; \mathbb{R}^{n\times d}).$ Put

$$ Y(u_k,t)=lim_{i\rightarrow +\infty} Y^i(u_k, t),\qquad  Z(u
_k,t)=lim_{i\rightarrow +\infty} Z^i(u_k, t).$$
Using similar steps as in the evaluation of the last limeses, and by applying  in addition Burkholder-Davis-Gundy inequality, one can obtain that
\begin{eqnarray*}
&& E\sup_{t\in[0,T]}| Y^i(u_k, t)-Y(u_k,t)|^2\rightarrow 0, \   {i\rightarrow +\infty}, \qquad \\
&& E\left(\int_0^T||Z^i(u_k, t)- Z(u_k,t)||^2ds\right)\rightarrow 0, \   {i\rightarrow +\infty}.
\end{eqnarray*}
So,  $\{(Y(u_k,t), Z(u_k,t)), t\in[0,T]\}$ solves Eq.(\ref{heavy1}) for final conditions $Y(u_k,T)=\xi(u_k),  k=\overline{1,N}$.\ \hskip 1.8cm $\star$
\end{dproof}

\smallskip

\underline{\it {Step 2.}} {\bf Solution to BSDE with interaction with discrete measure for all particles.}

\medskip
It is left to complete the result of existence and uniqueness of the solution of Eq.(\ref{heavy1}) for arbitrary $u\in \mathbb{R}^d$

\begin{equation}
Y(u,t)=\xi( u)-\int_t^T\sum_{j=1}^N h(Y(u,s),Y(u_j,s)) p_j ds-\int_t^TZ(u,s) dB(s), \label{anyfinal}
\end{equation}
where $Y(u_j,t), t\in [0,T], j=\overline{1,N}$ are solutions to the system of equations with discrete measure for heavy particles. Eq. (\ref{anyfinal}) has a drift which is dependent only of state process $Y(u,s)$ , and not  $Z(u,s)$, for $s\in[t,T]$, and satisfies Lipshitz condition. It follows from martignale representation theorem  that there exists its unique solution.
Also, as Eq. (\ref{anyfinal}) is as a type of Eq. (\ref{1}), there exists its unique solution by known result of Pardoux and Peng.

\medskip

\section{ Solution to BSDE with interaction for general measure}

\medskip

Aim of this section is to prove the existence and uniqueness of solution to

\begin{eqnarray}
&& d Y(u,t)=\int_{\mathbb{R}^n} h(Y,v) \mu_t(dv)dt+ Z(u,t) dB(t),\nonumber \\
&&  Y(u,T)=\xi(u),\nonumber\\
&&\mu_t=\mu_0\circ Y^{-1}(\cdot,t), \label{heavy_par1}
\end{eqnarray}
for general measure $\mu_0$. For this we will use approximation of arbitrary  measure by discrete measures 

$$
\sum_{k=1}^N p_k \delta_{u_k}\rightarrow \mu_0,
$$
in Wasserstein distance. Set $C(\mu, \mu_1)$ for measures  $\mu,\mu_1$ on $ \mathbb{R}^d$  ($\mu,\mu_1$  are from  \textfrak{M})is the set of all probability measures on the Borel $\sigma$-field in $  \mathbb{R}^{2d}$  which have $\mu$ and $\mu_1$ as its marginal projections.



\begin{definition}\label{wd0}
The Wasserstein distance of order zero between measures $\mu$ and $\mu_1$  is the value 
$$\gamma_0(\mu, \nu)=\inf_{  \varkappa\in C(\mu, \mu_1)} \int _{\mathbb{R}^d}\int _{\mathbb{R}^d}\frac{|u-v|}{1+|u-v|}\varkappa(du,dv).$$
\end{definition}
For $n\geq 1$  let

$$\textfrak{M}_n=\{\mu \in \textfrak{M}|\ \  \forall u \in \mathbb{R}: \int_{\mathbb{R}^d}|u-v|^n \mu(dv)<+\infty\}. $$

\begin{definition}\label{wdn}
	The Wasserstein distance of order $n$ between measures $\mu$ and $\mu_1$  is the value 
	$$\gamma_n(\mu, \nu)=\Big(\inf_{  \varkappa\in C(\mu, \mu_1)} \int _{\mathbb{R}^d}\int _{\mathbb{R}^d} |u-v|^n\varkappa(du,dv)\Big)^ {\frac{1}{n}}.$$
\end{definition}

\begin{remark} It can be proven that $\gamma_0$ is metric, and $(\textfrak{M},\gamma_0)$ is a complete separable metric space. Furthermore, it can be shown that  $(\textfrak{M}_n,\gamma_n)$ is a complete separable metric space for all $n\geq 1$ (see Chapter 1 in \cite{du89}).

\end{remark}

\medskip

We need following result.

\begin{lemma} \label{lemah} Suppose that $h$ satisfies  assumption $(A2)$. Then for 

$$f(y,\mu)= \int_{\mathbb{R}^d} h( y,u) \mu(du),$$
there exist constant $C>0$, such that for every $y, y_1\in\mathbb{R}^d, \mu,\mu^1\in \textfrak{M}_2$ following holds


$$  	|f(y, \mu)-	f(y_1, \mu^1)|^2\leq C\big[|y-y_1| ^2+ \gamma_2^2(\mu,
\mu^1)\big].$$

\end{lemma}

\begin{dproof} Consider for $ \varkappa\in C(\mu, \mu_1)$ and applying Cauchy--Schwarz inequality, it follows that
\begin{eqnarray}
&&| f(y, \mu)-f(y_1, \mu^1)|^2\\
&&\qquad = \left|\int_{\mathbb{R}^d} h(y,v) \mu(dv)-\int_{\mathbb{R}^d} h(y_1,v) \mu(dv)+\int_{\mathbb{R}^d} h(y_1,v) \mu(dv)-\int_{\mathbb{R}^d} h( y_1,v) \mu^1(dv)\right|^2\nonumber\\
&&\qquad\leq 2|y-y_1|^2+ 2\left|\int_{\mathbb{R}^d}\int_{\mathbb{R}^d}\left |h( y_1,u) - h( y_1,v)\right| \varkappa(du,dv)\right|^2\nonumber\\
&&\qquad\leq 2|y-y_1|^2+ 2L^2\int_{\mathbb{R}^d}\int_{\mathbb{R}^d }|u-v|^2\varkappa(du,dv). \label{gamma_new1}
\end{eqnarray}
By applying $\inf$ with respect to  $ \varkappa$, one can conclude
$$|f(y, \mu)-f(y_1, \mu^1)|^2\leq C\big[|y-y_1| ^2+ \gamma_2^2(\mu,
\mu^1)\big],
$$
where  $C$ is generic constant. This completes the proof. \hskip 7.5cm $\star$

\end{dproof}

\begin{remark}
Note that solution  to Eq. (\ref{heavy_par1}) continuously depend on the measure $\mu_0$ and the mapping $\xi$.
\end{remark}

\medskip

Indeed for  two arbitrary measures $\mu^i\in  \textfrak{M}, i=1,2 $, 
\begin{eqnarray}
&& d Y^i(u^i,t)=\int_{\mathbb{R}^n} h(Y^i,v^i) \mu_t^i(dv^i)dt+ Z^i(u^i,t) dB(t),\nonumber \\
&&  Y^i(u,T)=\xi(u^i), \nonumber\\
&& \mu_t^i=\mu_0\circ (Y^i)^{-1}(\cdot,t) \label{heavy_par1}
\end{eqnarray}

\begin{lemma} \label{l2.10} Under assumption  {\it {\bf (A2)} }
	\begin{eqnarray}
&& \hskip-1.3cm E|Y^1(u^1,t)-Y^2(u^2,t)|^2+E\int_t^T||Z^1(u^1,s)-Z^2(u^2,s)||^2ds\nonumber \\
&&\qquad \qquad \qquad \qquad \leq \tilde{C}\left( E[|u_1-u_2|^2+\int_t^T\gamma_2^2(\mu^1_t,\mu^2_t)dt\right),\label{lemma2.10}
	\end{eqnarray}	
	where $\tilde{C}$ is generic constant.

	\end{lemma}

\begin{dproof} 
Taking first  integral form of (\ref{heavy_par1}), we obtain 
\begin{eqnarray*}
&& Y^1(u^1,t)-Y^2(u^2,t)=u_1-u_2\\
&&+\int_t^T[f(Y^1(u^1,s), \mu_s^1)\\
&&-	f(Y^1(u^2,s), \mu_s^2)]ds-\int_t^T(Z^1(u^1,s)-Z^2(u^2,s)) dB(s).
\end{eqnarray*}

By It\^o formula  and  Lemma \ref{lemah},
\begin{eqnarray}
&&\hskip-1cm |Y^1(u^1,t)-Y^2(u^2,t)|^2=|\xi(u_1)-\xi(u_2)|^2\nonumber\\
&&+2\int_t^T|Y^1(u^1,s)-Y^2(u^2,s)|[f(Y^1(u^1,s), \mu_s^1)-	f(Y^1(u^2,s), \mu_s^2)]ds\nonumber\\
&&-2\int_t^T|Y^1(u^1,s)-Y^2(u^2,s)|[Z^1(u^1,s)-Z^1(u^2,s)]dB(s)\nonumber\\\
&&-\int_t^T||Z^1(u^1,s)-Z^2(u^2,s)||^2dt.
\end{eqnarray}
Hence
\begin{eqnarray}
&& E|Y^1(u^1,t)-Y^2(u^2,t)|^2+E\int_t^T||Z^1(u^1,s)-Z^2(u^2,s)||^2ds\nonumber\\\
&&\qquad \qquad \leq L |u_1-u_2|^2+E \int_t^T|Y^1(u^1,s)-Y^2(u^2,s)|^2ds\nonumber\\
&& \phantom{|Y^1(u^1,t)-}+E\int_t^T|f(Y^1(u^1,s), \mu_s^1)-	f(Y^1(u^2,s), \mu_s^2)|^2ds\nonumber\\\
&& \phantom{|Y^1(u^1,t)-}-2E\int_t^T|Y^1(u^1,s)-Y^2(u^2,s)|[Z^1(u^1,s)-Z^1(u^2,s)]dB(s).\label{ito_lema}
\end{eqnarray}

Then
\begin{eqnarray*}
&& E|Y^1(u^1,t)-Y^2(u^2,t)|^2\leq  L |u_1-u_2|^2\\
&&\phantom{ E|Y^1(u^1,t)-Y^2(u^2,t)|^2\leq }+(C+1)E \int_t^T|Y^1(u^1,s)-Y^2(u^2,s)|^2ds\\
&&\phantom{ E|Y^1(u^1,t)-Y^2(u^2,t)|^2\leq }+C E\int_t^T  \gamma_2^2(\mu_s^1,
\mu_s^2)ds.
\end{eqnarray*}
If we apply by iteration Gronwall-Bellman lemma (see \cite{bainov}) on last expression (and having in mind (\ref{ito_lema})), we obtain (\ref{lemma2.10}). \hskip 14.5cm $\star$

\end{dproof}
\vskip 0.5cm
\noindent It is known that arbitrary measure $\mu_0\in\textfrak{M}_2$ can be approximated by discrete measures in a sense of $\gamma_2$  (about possibilities of approximation see \cite{du89}).  We analyse equation with arbitrary measure $\mu_0$
\begin{eqnarray}
&& d Y(u,t)=\int_{\mathbb{R}^n} h(Y,v) \mu_t(dv)dt+ Z(u,t) dB(t),\nonumber\\
&&Y(u,T)=\xi(u),\nonumber\\
&&  \mu_t=\mu_0\circ Y^{-1}(\cdot,t) \label{heavy_par111}
\end{eqnarray}
Recall the equation for heavy particles
\begin{equation}\label{stability1}
	d Y(u_k,t)=\sum_{j=1}^{N_i} h(Y(u_k,t),Y^i(u_j,t)) p_j dt+ Z(u_k,t) dB(t),\ \ Y(u_k,T)=\xi(u_k).
	\end{equation}
Following theorem allows to prove general result, i.e.  to pass trough the limit from discrete measures to arbitrary measure.  

\begin{theorem} Let $({Y}^n(u,t),{Z}^n(u,t))_{0\leq t \leq T}\in \mathcal{M}^2( \mathbb{R}^d\times [0,T]; \mathbb{R}^d)\times \mathcal{M}^2(  \mathbb{R}^d\times [0,T]; \mathbb{R}^{d\times d}) $ be the sequence of the  solution of equations of type  (\ref{stability1}) for $\mu^n_T$.  Suppose that $({Y}(u,t),{Z}(u,t))_{0\leq t \leq T}\in \mathcal{M}^2( \mathbb{R}^d\times [0,T]; \mathbb{R}^d)\times \mathcal{M}^2(  \mathbb{R}^d\times [0,T]; \mathbb{R}^{d\times d}) $ solves the BSDE with interaction  (\ref{heavy_par111}) with arbitrary final measure $\mu_T$. If the sequence of BSDE with discrete measure  and BSDE with interaction  with arbitrary measure have the same final condition, furthermore if

\begin{equation}\gamma_2^2(\mu^n_T,\mu_T)\rightarrow 0, \qquad n \rightarrow \infty. \label{convgamma}
\end{equation}
in $\textfrak{M}_2$, then for every $T>0$ and $u \in \mathbb{R}^d$ 

\begin{eqnarray}
&&E|Y(u,t)-Y^n(u,t)|^2  \rightarrow 0,  n\rightarrow \infty,\label{yconv}\\
&&E\int_0^T||Z(u,t)-Z^n(u,t)||^2dt   \rightarrow 0,  n\rightarrow \infty.\label{zconv}
\end{eqnarray}

\end{theorem}



\begin{dproof}  From Defintion \ref{wd0} 

\begin{eqnarray}
&& E \gamma_2^2 (\mu_t^n,\mu_t)\leq  E\left(\int _{\mathbb{R}^d}\int _{\mathbb{R}^d}|Y(u,t)-Y^n(v,t)|^2\varkappa(du,dv)\right)\nonumber\\
&& \phantom{E \gamma_0 (\mu_t^n,\mu_t)^2}\leq  2E \left(\int _{\mathbb{R}^d}\int _{\mathbb{R}^d}|Y(v,t)-Y^n(v,t)|^2\varkappa(du,dv)\right)\nonumber\\
 && \phantom{E \gamma_0 (\mu_t^n,\mu_t)^2\leq}  
 +2E \left(\int _{\mathbb{R}^d}\int _{\mathbb{R}^d}|Y(u,t)-Y(v,t)|^2\varkappa(du,dv)\right)\nonumber\\
 && \phantom{E \gamma_0 (\mu_t^n,\mu_t)^2}= 2E \left(\int _{\mathbb{R}^d}(|Y(u,t)-Y^n(u,t)|^2)\mu(du)\right)\nonumber\\
 && \phantom{E \gamma_0 (\mu_t^n,\mu_t)^2\leq}  
 +2 \int _{\mathbb{R}^d}\int _{\mathbb{R}^d}
E |Y(u,t)-Y(v,t)|^2\varkappa(du,dv),\label{new1}
 \end{eqnarray}
 where in last equality we applied Fubini theorem (regarding that functions under integral are positive). For the estimate of first member from the right side of expression (\ref{new1}), we will use Lemma \ref{l2.10} in case we have the same final conditions (application of lemma with $u^1=u^2$), i.e.
 \begin{equation}
 2E \left(\int _{\mathbb{R}^d}(|Y(u,t)-Y^n(u,t)|^2)\mu(du)\right)\leq  \bar{C}E\int_t^T\gamma_2^2(\mu_t,\mu^n_t)dt,\label{1stmem}
 \end{equation}
 for some generic constant $ \bar{C}$.
 
For the second member we will  use the same lemma, but  in case when the measures are the same,
 $$E |Y(u,t)-Y(v,t)|^2\leq C' |u-v|^2,$$
 for some generic constant  $C'$. Then  substituting last estimates and (\ref{1stmem})  in (\ref{new1}), for $\varkappa\in C(\mu, \mu_1)$  one can get
 \begin{eqnarray*}
 &&E \gamma_2^2 (\mu_t^n,\mu_t)\leq  \bar{C}E\int_t^T\gamma_2^2(\mu_t,\mu^n_t)dt+ C'   \int _{\mathbb{R}^d}\int _{\mathbb{R}^d}E|u-v|^2\varkappa(du,dv) \\
&&\phantom{ E\gamma_0^2 (\mu_t^n,\mu_t)}\leq  \bar{C}'\left[ E\gamma_2^2(\mu_T,\mu^n_T)+ E\int_t^T\gamma_2^2(\mu_t,\mu^n_t)dt\right],
\end{eqnarray*}
where $ \bar{C}'=\max\{\bar{C}, C'\}$.
 Applying  Gronwall-Bellman's lemma to last inequality one can get

 \begin{equation}
  E\gamma_2^2 (\mu_t^n,\mu_t)\leq \bar{C}_1 T E\gamma_2^2(\mu_T^n,\mu_T),\label{gamma2}
  \end{equation}
where $\bar{C}_1$ is a generic constant.

From  Lemma \ref{l2.10} it follows that (\ref{yconv}) and (\ref{zconv}) hold. \hskip 6.5cm $\star$
 \end{dproof}
 
 \begin{theorem}
 There exist unique solution of BSDE with interaction (\ref{heavy_par111}) with arbitrary measure.
 \end{theorem}
 \begin{dproof} The proof of existence follows  from previous theorem, while uniqueness follows from Lemma  \ref{l2.10}. \hskip 13cm $\star$
  \end{dproof}

  \section{Conclusions and remarks}

  After giving the main results, we deduced following conclusions which are stated as remarks (for outgoing, as well as future work).

\begin{remark}
Eq. (\ref{1in}) can have a more general form with a function dependent of state process  $Y(\cdot)$ in diffusion coefficient additively added to process $Z(\cdot)$. Furthermore, diffusion coefficient can be generalized (as in the form of (\ref{add3})) to the case when it has some general function $b$ dependent of state process and measure $\mu$. Function $b$ should satisfy Lipschitz condition in the same manner as the drift coefficient.
\end{remark}

\newpage

{\bf Acknoledgements.} Jasmina \DJ or\dj evi\'c is supported by Grant STORM-Stochastics for Time-Space Risk Models, granted by Research Council of Norway - Independent projects: ToppForsk. Project nr. 274410.

\end{document}